\documentclass{amsart}

\newtheorem{theorem}{Theorem}[section]
\newtheorem{lemma}[theorem]{Lemma}

\theoremstyle{definition}
\newtheorem{definition}[theorem]{Definition}
\newtheorem{example}[theorem]{Example}

\newtheorem{corollary}[theorem]{Corollary}

\theoremstyle{remark}

\numberwithin{equation}{section}

\begin{document}

\title{Local Complete Intersections in $\mathbb{P}^2$ and Koszul 
syzygies}

\author{David Cox}
\address{Department of Mathematics and Computer Science, Amherst 
College, Amherst, MA 01002-5000}
\email{dac@cs.amherst.edu}
%\thanks{The first author was supported in part by NSF Grant 
%\#000000.}
\author{Hal Schenck}
\address{Department of Mathematics, Harvard University, Cambridge, MA 
02138; Current address: Department of Mathematics, Texas A\&M
University, College Station, TX, 77843}
\email{schenck@math.tamu.edu}
\thanks{The second author was supported by an NSF postdoctoral 
research fellowship.}

%    General info


\begin{abstract}
We study the syzygies of a codimension two ideal $I=\langle
f_1,f_2,f_3\rangle \subseteq k[x,y,z]$.  Our main result is that the
module of syzygies vanishing (scheme-theoretically) at the zero locus
$Z = {\bf V}(I)$ is generated by the Koszul syzygies iff $Z$ is a
local complete intersection.  The proof uses a characterization of
complete intersections due to Herzog \cite{h}.  When $I$ is saturated,
we relate our theorem to results of Weyman \cite{w} and Simis and
Vasconcelos \cite{sv}.  We conclude with an example of how our theorem
fails for four generated local complete intersections in $k[x,y,z]$
and we discuss generalizations to higher dimensions.
\end{abstract}

\maketitle

\section{Introduction}

Let $R = k[x,y,z]$ be the coordinate ring of $\mathbb{P}^{2}$, and
consider the ideal $I=\langle f_{1},f_{2},f_{3}\rangle \subseteq R$,
where $f_{i}$ is homogeneous of degree $d_{i}$.  It is well known that
the $f_{i}$ form a regular sequence in $R$ if and only if the Koszul
complex
\begin{equation}
\label{koszul}
\begin{aligned}
0 \longrightarrow R(-\sum\limits_{j=1}^3 d_j) &\xrightarrow{\left[ \!
\begin{array}{c}
f_3  \\
-f_2 \\
f_1 
\end{array}\! \right]} \bigoplus\limits_{j<k}R(-d_j-d_k)
\xrightarrow{\left[ \!\begin{array}{ccc}
f_2 & f_3&0 \\
-f_1 & 0 &f_3\\
0 & -f_1 & -f_2 
\end{array} \!\right]}\\ 
{}&\bigoplus\limits_{j=1}^3 R(-d_j)
\xrightarrow{\left[\! \begin{array}{ccc}
f_1 & f_2& f_3
\end{array} \!\right]} I \longrightarrow 0 
\end{aligned}
\end{equation}
is exact.  (As usual $R(i)$ denotes a shift in grading, i.e., $R(i)_j = 
R_{i+j}$.)

This paper will study the situation when $I = \langle
f_{1},f_{2},f_{3}\rangle$ has codimension two in $R$.  Thus $Z = {\bf
V}(I) \subseteq \mathbb{P}^{2}$ is a zero-dimensional subscheme. 
We will call $Z$ the \emph{base point locus} of $f_{1},f_{2},f_{3}$.

Since $Z \ne \emptyset$, $f_{1},f_{2},f_{3}$ no longer form a regular
sequence, so that \eqref{koszul} fails to be exact.  We can pinpoint
the failure of exactness as follows.

\begin{lemma}
\label{exactflop}
If $I = \langle f_{1},f_{2},f_{3}\rangle$ has codimension two in $R$,
then \eqref{koszul} is exact except at $\oplus_{j=1}^{3} R(-d_{j})$. 
In particular, the Koszul complex of $f_{1},f_{2},f_{3}$ gives the
exact sequences
\[
0 \to R(-\sum_{j=1}^3 d_j) \to
\bigoplus_{j<k}R(-d_j-d_k)
\to \bigoplus_{j=1}^3 R(-d_j)\ \ \text{and} \ \
\bigoplus_{j=1}^3 R(-d_j)
\rightarrow I \rightarrow 0.
\]
\end{lemma}

\begin{proof}
The exactness of the first sequence follows from $\mathrm{depth}(I) =
\mathrm{codim}(I) = 2$ and the Buchsbaum-Eisenbud exactness criterion
\cite{be}, and the second sequence is obviously exact.
\end{proof}

To describe how \eqref{koszul} behaves at 
$\bigoplus_{j=1}^{3} R(-d_{i})$, we make the following definition.

\begin{definition}
A \emph{Koszul syzygy} on $f_1,f_2,f_3$  is an element of the submodule
\[
K \subseteq \bigoplus_{j=1}^{3} R(-d_{i})
\]
generated by the columns of the matrix
\[
\left[\!\begin{array}{ccc} f_2 &
f_3&0 \\
-f_1 & 0 &f_3\\
0 & -f_1 & -f_2 
\end{array}\! \right]. 
\]
\end{definition}

This leads to the following corollary of Lemma~\ref{exactflop}.

\begin{corollary}
\label{kseq}
If $I = \langle f_1, f_2,f_3\rangle$ is a codimension two ideal, then 
\[
0 \longrightarrow  R(-\sum_{j=1}^3 d_j) \longrightarrow
\bigoplus_{j<k}R(-d_j-d_k) \longrightarrow K 
\longrightarrow 0
\]
is exact.
\end{corollary}

By Lemma~\ref{exactflop}, $K$ is a proper submodule of the syzygy
module $S$ defined by the exact sequence
\begin{equation}
\label{sdef}
0 \longrightarrow S \longrightarrow \bigoplus_{j=1}^{3} R(-d_{j}) 
\longrightarrow I \longrightarrow 0.  
\end{equation}
Here is a simple example of how $K$ and $S$ can differ.

\begin{example}
\label{ex1}
Consider the ideal $I=\langle xy,xz,yz\rangle $. Then the syzygy module
$S$ is generated by the columns of the matrix
\[
\phi =\left[ \begin{array}{ccc}
z & 0 \\
-y & y \\
0 & -x 
\end{array} \right]. 
\]
The generators of $S$ are clearly not Koszul since they have degree 
1.  
\end{example}

In this example, note that the generators of $S$ do not vanish
on the base point locus $Z = \{(1,0,0),(0,1,0),(0,0,1)\}$.  Since the
Koszul syzygies clearly vanish at the base points, we get another way
to see that the generators are not Koszul.

The observation that Koszul syzygies vanish at the basepoints leads
one to ask the question:\ \emph{Is every syzygy which vanishes at the
base points a Koszul syzygy?} Before answering this question, we need
to recall what ``vanishing'' means.

\begin{definition}
\label{vanishdef}
A syzygy $(a_{1},a_{2},a_{3})$ on the generators of $I = \langle
f_{1},f_{2},f_{3}\rangle$ \emph{vanishes at the basepoint locus} $Z =
{\bf V}(I)$ if $a_i \in I^{sat}$ for all $i$.  We will let $V$ denote
the module of syzygies which vanish at the basepoint locus $Z$.  (See
\cite[Def.\ 5.4]{c} for a sheaf-theoretic version of this definition.)
\end{definition}

In this terminology, we have $K \subseteq V$ since Koszul syzygies
vanish at the base points, and the above question can be rephrased
as:\ \emph{Is $K = V$?} In Example \ref{ex1}, one can compute that $K
= V$, but $K \ne V$ can also occur---see \cite[(5.4)]{c} for an
example.  This leads to the central question of this paper:\
\emph{When is $K = V$?}

To answer this question, we need the following definition.

\begin{definition}
An ideal is a \emph{local complete intersection} (\emph{lci}) if it is
locally generated by a regular sequence.  A codimension two lci in
$\mathbb{P}^2$ is \emph{curvilinear} if it is locally of the form
$\langle x,y^k\rangle$, where $x,y$ are local coordinates.
\end{definition}

In \cite{c}, toric methods were used to prove $K = V$ whenever $I =
\langle f_1,f_2,f_3\rangle \subseteq R$ is a codimension two
curvilinear local complete intersection.  The paper \cite{c} also
conjectured that $K = V$ should hold in the lci case.  Our main
result is the following strong form of this conjecture.

\begin{theorem}
\label{main}
If $I = \langle f_1,f_2,f_3\rangle \subseteq R$ has codimension 
two, then $K = V$ if and only if $I$ is a local complete intersection.
\end{theorem}

The proof of this theorem will be given in the next section. 

\section{Proof of the Main Theorem }

We begin with some preliminary definitions and results which will be
used in the proof of Theorem \ref{main}.  As in Section 1, $R =
k[x,y,z]$ with the usual grading.

\begin{definition} A submodule $M$ of a finitely generated graded free
$R$-module $F$ is \emph{saturated} if
\[
M = \{u \in F \mid \langle x,y,z\rangle u \subseteq M\}.
\]
\end{definition}

Recall that the \emph{Hilbert polynomial} $H(M)$ of a finitely
generated graded $R$-module $M$ is the unique polynomial such
that 
\[
H(M)(n) = \dim_{k} M_{n}
\]
for $n \gg 0$, where $M_{n}$ is the graded piece of $M$ in degree $n$. 

We omit the straightforward proof of the following result.

\begin{lemma}
\label{sateq}
Let $M \subseteq N$ be saturated submodules of a finitely generated 
graded free $R$-module.  Then $M = N$ if and only if $H(M) = H(N)$.
\end{lemma}

Given a codimension two ideal $I = \langle f_{1},f_{2},f_{3}\rangle 
\subseteq R$, we define the modules $K$ and $V$ as in the previous 
section.  These modules satisfy
\[
K \subseteq V \subseteq \bigoplus_{j=1}^{3} R(-d_{j}).
\]
We next show that they are saturated.

\begin{lemma}
\label{kvsat}
$K$ and $V$ are saturated submodules of $\bigoplus_{j=1}^{3}
R(-d_{j})$.
\end{lemma}

\begin{proof}
We first consider $V$.  Definition~\ref{vanishdef} implies that
\begin{equation}
\label{vformula}
V = S\cap \Big(\bigoplus_{j=1}^{3} I^{sat}(-d_{j})\Big).
\end{equation}
It is easy to see that $S$ is saturated.  Since the intersection of 
saturated submodules is saturated, we conclude that $V$ is saturated.

Turning to $K$, suppose that $u \in \bigoplus_{j=1}^{3} R(-d_{j})$
satisfies $\langle x,y,z\rangle u \subseteq K$.  We may assume that
$u$ is homogeneous.  We need to prove that $u \in K$.  For this 
purpose, let $L = K + Ru$.  Then consider the short exact sequence 
\[
0 \longrightarrow K \longrightarrow L \longrightarrow L/K 
\longrightarrow 0,
\]
and note that $L/K$ has finite length since $L_{n} = K_{n}$ for $n \ge
\deg(u) + 1$.  Let $\mathfrak{m} = \langle x,y,z\rangle$ denote the
irrelevant ideal of $R$.  We obtain a long exact sequence in local
cohomology
\[
0 \longrightarrow H^0_{\mathfrak{m}}(K) \longrightarrow 
 H^0_{\mathfrak{m}}(L) \longrightarrow 
 H^0_{\mathfrak{m}}(L/K) \longrightarrow 
 H^1_{\mathfrak{m}}(K) \longrightarrow \cdots
\]
>From \cite[App.\ A.4]{e}, we know that for a graded $R$-module $M$ the
local cohomology $H^i_{\mathfrak{m}}(M)$ vanishes for $i <$ depth $M$,
so the exact sequence of Corollary \ref{kseq} implies that
$H^0_{\mathfrak{m}}(K) = H^1_{\mathfrak{m}}(K) = \{0\}$.  Since $L
\hookrightarrow \bigoplus_{j=1}^{3} R(-d_{j})$,
$H^0_{\mathfrak{m}}(L)=\{0\}$.  This forces
\[
H^0_{\mathfrak{m}}(L/K,R) = \{0\}.
\]
However, since $L/K$ is a module 
of finite length, we also have $H^0_{\mathfrak{m}}(L/K)=L/K$.  We
conclude that $L/K = \{0\}$, so that $K = L$.  This implies $u \in K$,
as desired.
\end{proof}

The final ingredient we need for the proof is the following result of 
Herzog \cite{h} which characterizes complete intersections in the 
local case.

\begin{theorem}
\label{herzog}
Let $\mathcal{O}_{p}$ be the local ring of a point $p \in
\mathbb{P}^{2}$, and let $\mathcal{I}_{p} \subseteq \mathcal{O}_{p}$
be a codimension two ideal.  Then
\begin{equation}
\label{herzogineq}
\dim_k\mathcal{I}_p/\mathcal{I}_p^2 \ge 2 \dim_k 
\mathcal{O}_p/\mathcal{I}_p.
\end{equation}
Furthermore, equality holds if and only if $\mathcal{I}_{p}$ is a
complete intersection in $\mathcal{O}_{p}$.
\end{theorem}

We can now prove our main result.

\begin{proof}[Proof of Theorem \ref{main}]
By Lemmas \ref{sateq} and \ref{kvsat}, we know that $K = V$ if and 
only if they have the same Hilbert polynomials.  In other words,
\begin{equation}
\label{keqv}
K = V \iff H(K) = H(V).
\end{equation}
We next compute $H(K)$ and $H(V)$.
\medskip

By Lemma \ref{kseq}, $H(K)$ is given by
\[
H(K) = \sum\limits_{j<k}H(R(-d_j-d_k))-H(R(-{\textstyle 
\sum_{j=1}^3} d_{j})).
\]
Using the exactness of the Koszul complex of the regular sequence
$x^{d_{1}},y^{d_{2}},z^{d_{3}}$, this formula simplifies to
\begin{equation}
\label{hk}
H(K) = \sum_{j=1}^{3}H(R(-d_j))-H(R).
\end{equation}

We next consider $H(V)$. Using \eqref{vformula} and the exact sequence
\eqref{sdef}, we obtain the exact sequence
\[
0 \longrightarrow V \longrightarrow 
\bigoplus\limits_{j=1}^3I^{sat}(-d_j) 
\longrightarrow I \cdot I^{sat} \longrightarrow 0.
\]
This gives the Hilbert polynomial
\[
H(V) = \sum_{j=1}^3H(I^{sat}(-d_j))-H(I \cdot I^{sat}).
\]
However, since $0 \to I^{sat} \to R \to R/I^{sat} \to 0$ is exact and
${\bf V}(I^{sat}) = {\bf V}(I) = Z$ is zero dimensional, we also have
\[
H(I^{sat}) = H(R) - \deg Z.
\]
Hence the above formula for $H(V)$ can be written
\[
H(V) = \sum_{j=1}^3H(R(-d_j))-3\deg Z - H(I \cdot I^{sat}).
\]

We next observe that $I^{2}$ and $I \cdot I^{sat}$ have the same
saturation.  To prove this, note that $I^2 \subseteq I \cdot I^{sat}$,
so that it suffices to show that $I \cdot I^{sat} \subseteq
(I^2)^{sat}$.  Let $f \in I \cdot I^{sat}$, so $f = \sum_{i=1}^k f_i
g_i$ with $f_i \in I, g_i \in I^{sat}$.  But then there exists an $n$
such that $\langle x,y,z\rangle^n g_i \subseteq I$ for all $i$, hence
$f \in (I^2)^{sat}$.

It follows that $H(I^2) = H(I \cdot I^{sat})$ since the two ideals
have the same saturation.  This allows us to write the above formula
for $H(V)$ in the form
\[
H(V) = \sum_{j=1}^3H(R(-d_j))-3\deg Z - H(I^{2}).
\]
The exact sequences $0 \to I^{2} \to R \to R/I^{2} \to 0$ and $0 \to
I/I^{2} \to R/I^{2} \to R/I \to 0$ show that
\[
H(I^{2}) = H(R) - H(R/I^{2}) = H(R) - H(I/I^{2}) - \deg Z.
\]
Then the previous formula for $H(V)$ can be written as
\begin{equation}
\label{hv}
H(V) = \sum_{j=1}^3H(R(-d_j))-H(R) + H(I/I^{2})-2\deg Z.
\end{equation}
Comparing this to \eqref{hk}, we obtain
\begin{equation}
\label{hkhv}
H(K) = H(V) \iff H(I/I^{2}) = 2\deg Z.
\end{equation}

It remains to compute $\deg Z$ and $H(I/I^{2})$. If $\mathcal{I}$ is
the ideal sheaf of $Z$, then
\[
\deg Z = \dim_k H^0(Z, \mathcal{O}_{Z}) = \dim_k H^0(\mathbb{P}^2,
\mathcal{O}_{\mathbb{P}^{2}}/\mathcal{I}) = \sum_{p \in Z} \dim_k
\mathcal{O}_p/\mathcal{I}_p,
\]
where $\mathcal{O}_p = \mathcal{O}_{\mathbb{P}^{2},p}$ and
$\mathcal{I}_p$ is the localization of $I$ at $p \in Z$.  Since
$I/I^2$ has zero dimensional support, we also have
\[
H(I/I^2) = \dim_k H^0(\mathbb{P}^2,\mathcal{I}/\mathcal{I}^2) = 
\sum_{p \in Z}\dim_k \mathcal{I}_p/\mathcal{I}_p^2.
\]
By \eqref{herzogineq}, we know that
\[
\dim_k \mathcal{I}_p/\mathcal{I}_p^2 \ge 2\dim_k 
\mathcal{O}_p/\mathcal{I}_p 
\]
for every $p \in Z$.  It follows easily that
\[
H(I/I^{2}) = 2\deg Z \iff 
\dim_k \mathcal{I}_p/\mathcal{I}_p^2 = 2\dim_k 
\mathcal{O}_p/\mathcal{I}_p\ \text{for all}\ p \in Z,
\]
and by the final assertion of Theorem \ref{herzog}, we conclude that
\begin{equation}
\label{lcicond}
H(I/I^{2}) = 2\deg Z \iff I\ \text{is lci}.
\end{equation}
Theorem \ref{main} follows immediately from \eqref{keqv}, 
\eqref{hkhv} and \eqref{lcicond}.
\end{proof}

\section{The Saturated Case}

This section will consider the case where $I \subset R = k[x,y,z]$ is
a saturated ideal of codimension two. Our goal is to explain how $I$
being a local complete intersection relates to the results of Weyman
\cite{w} and Simis and Vasconcelos \cite{sv}.

With the above hypothesis, $I$ is Cohen-Macaulay (see \cite[Proposition
5.2]{c}) of codimension $2$, so that by the Hilbert-Burch Theorem (see
\cite{e}), the free resolution of $I$ has the form
\[
0\longrightarrow \bigoplus_{i=1}^{m}R(-b_i)  
\xrightarrow{\ \phi\ } \bigoplus_{j=1}^{m+1}R(-a_j) 
\xrightarrow{\ \psi\ } I \longrightarrow 0.
\]
We will write this resolution as
\begin{equation}
\label{fgi}
0\longrightarrow F \xrightarrow{\ \phi\ } G \xrightarrow{\ \psi\ } I 
\longrightarrow 0,
\end{equation}
where $F$ and $G$ are free graded $R$-modules of ranks $m$ and $m+1$
respectively.  

We first consider the results of \cite{w}.  From \eqref{fgi}, Weyman
constructs the complex
\begin{equation}
\label{ressym2}
0 \longrightarrow \wedge^2 F \longrightarrow F \otimes G
\longrightarrow \mathrm{Sym}_2\, G \longrightarrow \mathrm{Sym}_2\, I
\longrightarrow 0.
\end{equation}
By Theorem 1 of \cite{w}, it follows easily that this complex is exact
in our situation.

Using the natural map $\mathrm{Sym}_2\, I \to I^{2}$, \eqref{ressym2} 
gives the complex
\begin{equation}
\label{resi2}
0 \longrightarrow \wedge^2 F \longrightarrow F \otimes G 
\longrightarrow \mathrm{Sym}_2\, G \longrightarrow I^{2} 
\longrightarrow 0.
\end{equation}
The results of Simis and Vasconcelos \cite{sv} lead to the following
theorem.

\begin{theorem}
\label{satthm}
Let $I \subseteq R$ be Cohen-Macaulay of codimension two.  Then the 
following conditions are equivalent:
\begin{enumerate}
    \item $I$ is a local complete intersection.
    \item The natural map $\mathrm{Sym}_{2}\,I \to I^{2}$ is an
    isomorphism. 
    \item $\mathrm{Sym}_{2}\, I$ is torsion free.
    \item The complex \eqref{resi2} is a free resolution of $I^{2}$.
\end{enumerate}
\end{theorem}

\begin{proof}
The paper \cite{sv} defines an invariant $\delta(I)$ and, in Corollary
1.2, shows that $\delta(I) = \ker(\mathrm{Sym}_{2}\,I \to I^2)$.
The first remark following Theorem 2.2 in \cite{sv} implies that
$\delta(I) = \{0\}$ if and only if $I$ is generically a complete
intersection (i.e., the localization of $R/I$ at each minimal prime is
a complete intersection).  Since $\mathbf{V}(I)$ has dimension zero,
the latter is equivalent to being a local complete intersection.  This
proves $(1) \iff (2)$.

The equivalence $(2) \iff (3)$ is obvious since $\mathrm{Sym}_2\, I
\to I^{2}$ is onto and its kernel is torsion by \cite[Ex.\ A2.4, p.\
574]{e}.  Finally, $(2) \iff (4)$ follows from the exactness of
\eqref{ressym2}.
\end{proof}

We conclude this section by showing that there are geometrically
interesting ideals $I \subset R$ which are codimension two
Cohen-Macaulay local complete intersections.  The ideals we will
consider arise in the study of line arrangements in $\mathbb{P}^2$.
Let
\[
Q = \prod_{i=1}^d \ell_i,
\]
where the $\ell_i$ are distinct linear forms in $R = k[x,y,z]$. 
Localizing and applying the Euler relation shows that the Jacobian ideal 
\[
J_Q = \langle Q_{x},Q_{y},Q_{z}\rangle \subseteq R
\]
is a local complete intersection, generated by three forms of degree
$d-1$.  An open problem in the study of line arrangements
consists in determining when $R/J_Q$ is Cohen-Macaulay.  In
this situation the arrangement is called a \emph{free} arrangement;
for more on this, see \cite{ot}.

There are certain cases where freeness is known.  For 
example, let 
\begin{align*}
L_1 &=\prod_{i=1}^m y-a_i x\\
L_2 &=\prod_{j=1}^n z-b_j x,
\end{align*}
where $a_i, b_j$ are nonzero and distinct.  Then put
\[
Q =x L_1 L_2.
\]
The Addition Theorem (\cite[Thm.\ 4.50]{ot}) implies that these
arrangements are free.  It follows that the corresponding Jacobian
ideal $J_{Q}$ is a Cohen-Macaulay local complete intersection, and it has
codimension two since ${\bf V}(J_{Q})$ is the singular locus of ${\bf
V}(Q)$, which is union of distinct lines.

For these arrangements, the free resolution is given by 
\[
0\longrightarrow R(-m-2n) \oplus R(-2m-n)
\longrightarrow R^3(-m-n) \longrightarrow J_Q \longrightarrow 0,
\]
One can show that $H(R/J_Q) = m^2+n^2+mn$ and 
$H(J_{Q}/J_Q^2) = 2(m^2+n^2+mn)$.  Notice how these numbers are 
consistent with the proof of Theorem~\ref{main}.  Also observe that 
Example~\ref{ex1} from Section 1 is the special case of this 
construction corresponding to $L_{1} = y$ and $L_{2} = z$.

\section{Final Remarks}

We first give an example to show that Theorem \ref{main} can fail if
an ideal in $R = k[x,y,z]$ has more than three generators.

\begin{example}
Let $J \subseteq R$ be the ideal of all forms vanishing on five
general points in $\mathbb{P}^2$.  Then $J$ is generated by a conic
and two cubics, so the degree three piece of $J$ has five generators. 
If we take four generic elements $f_{1},f_{2},f_{3},f_{4} \in J_3$,
then one can show the following:
\begin{itemize}
\item $I = \langle f_{1},f_{2},f_{3},f_{4}\rangle$ and $J$ define
the same subscheme of $\mathbb{P}^2$.  It follows that $I$ is an lci of 
codimension two.
\item Not every syzygy on $f_{1},f_{2},f_{3},f_{4}$ which vanishes
at the basepoints is Koszul.
\end{itemize}
\end{example}

So three generators is the number needed to make things work for 
codimension two ideals of $R = k[x,y,z]$.  

In proving Theorem~\ref{main}, we showed that the Hilbert polynomial
of the Koszul syzygies $K$ is determined by the degrees of the
generators (see \eqref{hk}), while we need more subtle information in
order to compute the Hilbert polynomial of the syzygies $V$ vanishing
at the basepoints (see \eqref{hv}).  The miracle is that if $I$ is a
three generated ideal of codimension two, then exactly the right
numerology occurs so that $V = K$ is equivalent to being a local
complete intersection.

It makes sense to ask if a similar ``numerology'' occurs in higher
dimensions.  A natural place to start would be to consider ideals
\[
I = \langle f_{1},\dots,f_{n+1}\rangle \subset k[x_{0},\dots,x_{n}]
\]
whose base point locus has codimension $n$.  In the local case, such
ideals are called \emph{almost complete intersections}.  It would be
interesting to see if the results of Sections 2 and 3 extend to this
situation.

Finally, we should mention that evidence for Theorem \ref{main} was
provided by many Macaulay2 computations. Macaulay2 is available at 
the URL
\[
{\tt http://www.math.uiuc.edu/Macaulay2/}
\]
We are also grateful to Winfried Bruns, J\"urgen Herzog and Aron Simis
for helpful comments.

\bibliographystyle{amsplain}

\end{document}